\pdfoutput=1
\RequirePackage{ifpdf}
\ifpdf % We are running pdfTeX in pdf mode
\documentclass[pdftex]{sigma}
\else
\documentclass{sigma}
\fi

\usepackage{bm}

\newtheorem{Theorem}{Theorem}[section]
\newtheorem{Conjecture}[Theorem]{Conjecture}
\newtheorem{Lemma}[Theorem]{Lemma}
\numberwithin{equation}{section}

\renewcommand{\phi}{\varphi}
\newcommand{\C}{\operatorname{C}}

\newcommand{\ZZ}{\mathbb{Z}}
\newcommand{\CC}{\mathbb{C}}

\newcommand{\NN}{\mathbb{N}}
\newcommand{\FF}{\mathbb{F}}
\newcommand{\MP}{\mathcal{P}}
\newcommand{\GL}{\operatorname{GL}}
\newcommand{\Irr}{\operatorname{Irr}}
\newcommand{\IBr}{\operatorname{IBr}}
\newcommand{\sgn}{\operatorname{sgn}}
\newcommand{\diag}{\operatorname{diag}}

\begin{document}

\renewcommand{\thefootnote}{}

\renewcommand{\PaperNumber}{100}

\FirstPageHeading

\ShortArticleName{Morita Equivalent Blocks of Symmetric Groups}

\ArticleName{Morita Equivalent Blocks of Symmetric Groups\footnote{This paper is a~contribution to the Special Issue on the Representation Theory of the Symmetric Groups and Related Topics. The full collection is available at \href{https://www.emis.de/journals/SIGMA/symmetric-groups-2018.html}{https://www.emis.de/journals/SIGMA/symmetric-groups-2018.html}}}

\Author{Benjamin SAMBALE}

\AuthorNameForHeading{B.~Sambale}

\Address{Fachbereich Mathematik, TU Kaiserslautern, 67653 Kaiserslautern, Germany}
\Email{\href{mailto:sambale@mathematik.uni-kl.de}{sambale@mathematik.uni-kl.de}}

\ArticleDates{Received April 16, 2018, in final form September 18, 2018; Published online September 20, 2018}

\Abstract{A well-known result of Scopes states that there are only finitely many Morita equivalence classes of $p$-blocks of symmetric groups with a given weight (or defect). In this note we investigate a lower bound on the number of those Morita equivalence classes.}


\Classification{20C30; 20C08}

\renewcommand{\thefootnote}{\arabic{footnote}}
\setcounter{footnote}{0}

\section{Introduction}

Every $p$-block $B$ of a symmetric group $S_n$ is uniquely determined by its \emph{weight} $w\ge 0$ and its \emph{core} $\mu$ (a partition of $n-pw$, see next Section~\ref{section3} for details). Scopes~\cite{Scopes} showed that there are only finitely many Morita equivalence classes of $p$-blocks of symmetric groups over the field $\FF_p$ with a~given weight (this confirms a special case of Donovan's conjecture). Puig~\cite{PuigScopes} extended her result to blocks defined over~$\ZZ$. More precisely, the number of these so-called \emph{Scopes classes} is $\frac{1}{p}\binom{wp}{p-1}$ provided $w>0$. If $p>2$ and $w>1$, then there exist pairs of distinct Scopes classes whose cores are conjugate to one another. Since the corresponding blocks of such pairs are Morita equivalent (even isomorphic), the number of Morita equivalence classes $M(p,w)$ of $p$-blocks with weight $w$ satisfies
\begin{gather}\label{bound}
M(p,w)\le\frac{1}{2p}\binom{wp}{p-1}+\frac{1}{2}\binom{\lfloor wp/2\rfloor}{\lfloor p/2\rfloor}
\end{gather}
for $w>0$ (see \cite[Corollary~3.10]{Richards}). In this paper we are interested in the sharpness of this bound. Obviously, there is only one class of blocks with weight (defect)~$0$ or~$1$ respectively. The blocks of weight $2$ have been investigated by many authors. For instance, Richards~\cite{Richards} obtained a formula for the decomposition matrix of those blocks. Since the decomposition matrix is a~Morita invariant over $\ZZ$ (up to permutations of rows and columns), it is possible to count Morita equivalence classes at least when $p$ is small. We have used the computer algebra system GAP~\cite{GAP48} to confirm the following conjecture for $p\le 11$.

\begin{Conjecture}\label{con1}
\begin{gather*}M(p,2)=\frac{1}{2p}\binom{2p}{p-1}+\frac{1}{2}\binom{p}{\lfloor p/2\rfloor}.\end{gather*}
\end{Conjecture}

For blocks of weight $w=3$ and $p\ge 5$, Fayers~\cite{Fayers3} has shown that all decomposition numbers are $0$ or $1$. Using the Jantzen--Schaper formula~\cite{Schaper} it is therefore possible to compute the decomposition matrices recursively. In this way we obtain $M(5,3)=147$ and $M(7,3)=3936$ using GAP. Note that \eqref{bound} is sharp in these cases as well.

Much less is known about the blocks of weight $w>3$. Nevertheless, for small primes the following complementary result seems to hold.

\begin{Conjecture}\label{main}
\begin{gather*}M(p,w)=\begin{cases}
2&\text{if } p=2,\ w=3,\\
w&\text{if } p=2,\ w>3,\\
\left\lfloor\dfrac{3w^2+2w}{4}\right\rfloor&\text{if } p=3,\ w>0.
\end{cases}\end{gather*}
\end{Conjecture}

Using an algorithm which is explained in the next sections, we have checked Conjecture~\ref{main} for $p=2$, $w\in\{0,\ldots,28\}\cup\{32\}$ and $p=3$, $w\le 15$. The case $(p,w)=(2,3)$ is in fact an exception. Here, the Scopes classes are represented by the cores $()$, $(1)$ and $(2,1)$. One can show with Magma~\cite{Magma} that the cores $(1)$ and $(2,1)$ correspond to Morita equivalent blocks at least over $\FF_2$.

Unfortunately, our method does not extend directly to $p>3$. Nevertheless, we obtain the following intervals with GAP: $496\le M(5,4)\le 507$ and $1278\le M(5,5)\le 1298$.

Finally, we remark that our observations should also apply to Hecke algebras.

\section{A Morita invariant}\label{section2}

Since the decomposition matrix of a block is usually not available (like in the situation of Conjecture~\ref{main}), it is important to know other Morita invariants. The following approach applies to any finite group~$G$ and any prime~$p$. We denote the set of $p$-regular elements of $G$ by $G^0$ and the set of irreducible Brauer characters by $\IBr(G)$. For $\chi,\psi\in\Irr(G)\cup\IBr(G)$ we are interested in the following $p$-scalar product
\begin{gather*}[\chi,\psi]^0:=\frac{1}{|G|}\sum_{g\in G^0}\chi(g)\psi\big(g^{-1}\big)\in\CC.\end{gather*}
The following lemma states that these numbers are preserved under Morita equivalence.

\begin{Lemma}\label{lem} Let $B_1$ and $B_2$ be Morita equivalent blocks of $($possibly different$)$ finite groups with respect to a complete discrete valuation ring. Let
\begin{gather*}M_i:=\big([\chi,\psi]^0\big)_{\chi,\psi\in\Irr(B_i)}\end{gather*}
for $i=1,2$. Then there exists a permutation matrix $T$ such that
\begin{gather*}TM_1=M_2T.\end{gather*}
\end{Lemma}
\begin{proof}Let $Q_i=(d_{\chi\psi})$ be the decomposition matrix of $B_i$ for $i\in\{1,2\}$. Then $C_i:=Q_i^\text{t}Q_i$ is the Cartan matrix of $B_i$ where $Q_i^\text{t}$ denotes the transpose of $Q_i$. By \cite[Theorem~2.13]{Navarro} we have $C_i^{-1}=\big([\phi,\theta]^0\big)_{\phi,\theta\in\IBr(B_i)}$. It follows that
\begin{gather}
M_i=\bigg(\frac{1}{|G|}\sum_{g\in G_i^0}\sum_{\phi,\theta\in\IBr(B_i)}d_{\chi\phi}d_{\psi\theta}\phi(g)\theta\big(g^{-1}\big)\bigg)_{\chi,\psi\in\Irr(B_i)}\nonumber\\
\hphantom{M_i}{} =\bigg(\sum_{\phi,\theta\in\IBr(B_i)}d_{\chi\phi}d_{\psi\theta}[\phi,\theta]^0\bigg)_{\chi,\psi\in\Irr(B_i)}=Q_iC_i^{-1}Q_i^\text{t}.\label{form}
\end{gather}
Since $B_1$ and $B_2$ are Morita equivalent, there exist permutation matrices $S$, $T$ such that \linebreak \smash{$TQ_1S=Q_2$}. Since $T$ is orthogonal, we conclude that
\begin{gather*}M_2T=Q_2\big(Q_2^\text{t}Q_2\big)^{-1}Q_2^\text{t}T=TQ_1S\big(S^{\text{t}}Q_1^\text{t}T^\text{t}TQ_1S\big)^{-1}S^\text{t}Q_1^\text{t}T^\text{t}T=TQ_1C_1^{-1}Q_1^\text{t}=TM_1.\!\!\!\tag*{\qed}
\end{gather*}\renewcommand{\qed}{}
\end{proof}

The matrix $M:=M_1$ in Lemma~\ref{lem} encodes several important numerical invariants of \smash{$B:=B_1$}. Obviously, the size of $M$ is $k(B)\times k(B)$ where $k(B):=|\Irr(B)|$. By \eqref{form}, the rank of $M$ equals $l(B):=|\IBr(B)|$. Recall that the \emph{height} $h(\chi)\ge 0$ of $\chi\in\Irr(B)$ is defined by $\chi(1)_p=p^{a-d+h(\chi)}$ where $|G|_p=p^a$ and $d$ is the defect of $B$. If $h(\chi)=0$, then the $p$-adic valuation of $[\chi,\psi]^0$ equals $h(\psi)-d$ for every $\psi\in\Irr(B)$ (see \cite[Theorem~3.24]{Navarro}). Moreover, if $h(\chi)>0$, then the $p$-adic valuation of $[\chi,\chi]^0$ is larger than $-d$ (see \cite[Lemma~3.22(a)]{Navarro}).
In this way, $M$ encodes all character heights. Finally, one can show that $M$ determines the elementary divisors of the Cartan matrix of $B$. However, the next theorem implies that all these invariants do not suffice to distinguish blocks of symmetric groups.

\begin{Theorem}[{Chuang--Rouquier~\cite[Theorem~7.2]{ChuangRouquier}}]\label{CR}
Blocks $B_1$ and $B_2$ of symmetric groups with the same weight are splendidly derived equivalent. In particular there exists a \emph{signed} permutation matrix $T$ such that $TM_1=M_2T$ in the situation of Lemma~{\rm \ref{lem}}.
\end{Theorem}

We will see in the following that for blocks of symmetric groups the matrix $M$ can be computed from the character table of a local subgroup which only depends on the weight of the block (see Theorem~\ref{osima}). At the same time we determine the signs of~$T$ in Theorem~\ref{CR} without using the characters explicitly. To do so, we need to introduce a lot of notation.

\section{Notation}\label{section3}

We fix the following notation (details can be found in \cite{OlssonCombi}). A \emph{partition} of $n\in\NN_0$ is a non-increasing sequence of positive integers $\lambda=(\lambda_1,\ldots,\lambda_l)$ such that $|\lambda|:=\sum\limits_{i=1}^l\lambda_i=n$. The number $l(\lambda):=l$ is called the \emph{length} of $\lambda$. We allow the empty partition~$()$ of~$0$. The set of partitions of~$n$ is denoted by~$\MP(n)$.

It is well-known that the conjugacy classes of the symmetric group $S_n$ of degree $n$ consist of the elements with a common cycle structure. Hence, we can choose a set of representatives $\{s_\lambda\colon \lambda\in\MP(n)\}$ of those conjugacy classes. Similarly, the irreducible characters of $S_n$ are parametrized by $\MP(n)$ and we will write \begin{gather*}\Irr(S_n)= \{\chi_\lambda\colon \lambda\in\MP(n) \}.\end{gather*}
Let $p$ be a prime. Successively removing all hooks of length $p$ from (the Young diagram of) $\lambda\in\MP(n)$ yields the \emph{$p$-core} $\lambda_{(p)}$ of $\lambda$. The number of removed hooks is called the \emph{weight} $w$ of~$\lambda$. Observe that $n=|\lambda_{(p)}|+pw$. Moreover, the \emph{$p$-sign} of~$\lambda$ is defined by $\delta_p(\lambda)=(-1)^{\sum l_i}$ where the $l_i$ are the leg lengths of the removed hooks. By Nakayama's conjecture, characters $\chi_\lambda,\chi_\mu\in\Irr(S_n)$ lie in the same $p$-block if and only if~$\lambda_{(p)}=\mu_{(p)}$. Therefore, we may speak of the core and the weight of a block~$B$ of~$S_n$.

Next we define the $p$-quotient of $\lambda\in\MP(n)$. Let $l=l(\lambda)$ and choose $s\ge 0$ such that $l+s\equiv 0\pmod{p}$. Define
\begin{gather*}\beta:=(\lambda_1-1+l+s,\,\lambda_2-2+l+s,\ldots,\lambda_l+s,\,s-1,\,s-2,\ldots,0)\end{gather*}
(this is a \emph{$\beta$-set} for $\lambda$). For $i=0,\ldots,p-1$ let $\beta_i:=\{x\in\NN_0\colon px+i\in\beta\}$. Writing the elements of $\beta_i$ in decreasing order $\beta_i=\{b_1,\ldots,b_k\}$ gives a partition $P(\beta_i):=(b_1-k+1,b_2-k+2,\ldots,b_k)$ where we omit parts which are zero. Finally, the $p$-tuple of partitions
\begin{gather*}\bm{\lambda}:=\bigl(P(\beta_0),\ldots,P(\beta_{p-1})\bigr)\end{gather*}
is called the \emph{$p$-quotient} of $\lambda$. It turns out that $\sum\limits_{i=0}^{p-1}|P(\beta_i)|$ equals the weight of $\lambda$. In general, the set of $p$-tuples $\bm{\lambda}=(\bm{\lambda}_0,\ldots,\bm{\lambda}_{p-1})$ of partitions such that $\sum|\bm{\lambda}_i|=w$ is denoted by $\MP^p(w)$. The elements of $\MP^p(w)$ will always be denoted by bold Greek letters.

Now we relate the characters of the block $B$ above to the characters of the wreath product $G_w:=C_p\wr S_w$ where $C_p$ is a (fixed) group of order $p$. Let $\Irr(C_p)=\{\psi_0=1,\ldots,\psi_{p-1}\}$ and $\bm{\lambda}\in\MP^p(w)$. For $i=0,\ldots,p-1$, the linear character $\psi^{\otimes|\bm{\lambda}_i|}:=\psi_i\otimes\cdots\otimes\psi_i\in\Irr(C_p^{|\bm{\lambda}_i|})$
has a~unique extension to $C_p\wr S_{|\bm{\lambda}_i|}$ (still denoted by $\psi^{\otimes|\bm{\lambda}_i|}$) which acts trivially on $S_{|\bm{\lambda}_i|}$. Then we can define $\psi_{\bm{\lambda}_i}:=\psi^{\otimes|\bm{\lambda}_i|}\chi_{\bm{\lambda}_i}\in\Irr(C_p\wr S_{|\bm{\lambda}_i|})$ where $\chi_{\bm{\lambda}_i}$ is the inflation from $S_{|\bm{\lambda}_i|}$. Finally let
\begin{gather*}\psi_{\bm{\lambda}}:=(\psi_{\bm{\lambda}_0}\otimes\cdots\otimes\psi_{\bm{\lambda}_{p-1}})^{G_w}.\end{gather*}
The degree of this character is
\begin{gather}\label{deg}
\psi_{\bm{\lambda}}(1)=\binom{w}{|\bm{\lambda}_0|,\ldots,|\bm{\lambda}_{p-1}|}\chi_{\bm{\lambda}_0}(1)\cdots\chi_{\bm{\lambda}_{p-1}}(1).
\end{gather}
By the hook formula and \cite[Lemma~2.1]{OlssonMcKay}, the map $\Irr(B)\to\Irr(G_w)$, $\chi_\lambda\mapsto\psi_{\bm{\lambda}}$ is a height preserving bijection, i.e., $p^w(w!)_p\chi_\lambda(1)_p=(n!)_p\psi_{\bm{\lambda}}(1)_p$ for every $\chi_\lambda\in\Irr(B)$.

In order to obtain information on the decomposition matrix of $B$ we need to label the conjugacy classes of $G_w$.
Here we identify $C_p$ with $\ZZ/p\ZZ$ to simplify notation. For $(x_1\cdots x_w,\sigma)\in G_w$ (with $x_1,\ldots,x_w\in\ZZ/p\ZZ$, $\sigma\in S_w$) we define $\bm{\lambda}\in\MP^p(w)$ as follows: For every cycle $(a_1,\ldots,a_s)$ in $\sigma$ let $s\in\bm{\lambda}_{x_{a_1}+\dots+x_{a_s}}$. It turns out that two elements of $G_w$ are conjugate if and only if they yield the same $\bm{\lambda}$. Let \begin{gather*}\big\{g_{\bm{\lambda}}\colon \bm{\lambda}\in\MP^p(w)\big\}\end{gather*}
be a set of representatives for the conjugacy classes of $G_w$. The class sizes can be computed with the formula
\begin{gather*}|{\C}_{G_w}(g_{\bm{\lambda}})|=\prod_{i=0}^{p-1}p^{l(\bm{\lambda}_i)}\big|{\C}_{S_{|\bm{\lambda}_i|}}(s_{\bm{\lambda}_i})\big|.\end{gather*}

\section{Osima's result}\label{section4}

\begin{Theorem}[Osima~{\cite[Theorem~8]{OsimaSym}}]\label{osima}
With the notation above let $Q$ be the decomposition matrix of $B$ and let $D:=\diag(\delta_p(\lambda)\colon \chi_\lambda\in\Irr(B))$. Moreover, let
\begin{gather*}\Gamma:=\bigl\{\bm{\lambda}\in\MP^p(w)\colon \bm{\lambda}_0=()\bigr\}\qquad\text{and}\qquad X:=\bigl(\psi_{\bm{\lambda}}(g_{\bm{\mu}})\bigr)_{\bm{\lambda}\in\MP^p(w),\,\bm{\mu}\in\Gamma}.\end{gather*}
Then there exists $S\in\GL(l(B),\CC)$ $($depending on $B)$ such that
\begin{gather}\label{equosima}
Q=DXS,
\end{gather}
where the row of $\chi_\lambda$ corresponds to the row of $\psi_{\bm{\lambda}}$.
\end{Theorem}

We remark that \eqref{equosima} does not depend on the labeling of $\IBr(B)$, since any permutation of Brauer characters can be realized by $S$. As in the proof of Lemma~\ref{lem}, we have
\begin{gather*}M:=\bigl([\chi_\lambda,\chi_\mu]^0\bigr)_{\chi_\lambda,\chi_\mu\in\Irr(B)}=DX\big(X^\text{t}X\big)^{-1}X^\text{t}D=D\overline{X}\big(X^\text{t}\overline{X}\big)^{-1}X^\text{t}D,\end{gather*}
where $\overline{X}$ denotes the complex conjugate of $X$. The second orthogonality relation for $G_w$ implies that
\begin{gather}\label{XX}
X^\text{t}\overline{X}=\diag\bigl(|{\C}_{G_w}(g_{\bm{\mu}})|\colon \bm{\mu}\in\Gamma\bigr).
\end{gather}
In this way, $M$ can be computed from the character table of $G_w$ and some elementary combinatorics for $D$. Doing so, we are able to distinguish the Scopes classes for $p=3$ and $w\le 15$. For $p=5$ and $w\in\{4,5\}$, some of the Scopes classes yield the same $M$ (see last section). Hence, we obtain only a lower bound on $M(5,w)$ in these cases. For $p=2$ we refine our method in the next section in order to deal with larger weights.

\section[The case $p=2$ in Conjecture~\ref{main}]{The case $\boldsymbol{p=2}$ in Conjecture~\ref{main}}\label{section5}

For $p=2$ we notice first that $X$ is integral and \eqref{XX} simplifies to
\begin{gather*}X^\text{t}X=\diag\bigl(2^{l(\mu)}|{\C}_{S_w}(s_{\mu})|\colon \mu\in\MP(w)\bigr).\end{gather*}
We now define a ``small'' submatrix $M^0$ of $M$ which suffices to distinguish the $2$-blocks with weight $w\in\{1,\ldots,28\}\cup\{32\}$. We do not know if this approach works in general for $p=2$.
Let $\Irr_0(B):=\{\chi_\lambda\in\Irr(B)\colon h(\chi_\lambda)=0\}$. Since $M$ encodes character heights, it is clear that the permutation matrix $T$ in Lemma~\ref{lem} permutes $\Irr_0(B)$ (the characters of height $0$ correspond to the diagonal entries of $M$ with the lowest $p$-adic valuation, see remark after Lemma~\ref{lem}). We may therefore replace $M$ by $\bigl([\chi_\lambda,\chi_\mu]^0\bigr)_{\chi_\lambda,\chi_\mu\in\Irr_0(B)}$. Let $w=\sum\limits_{i\ge 0}a_i2^i$ be the $2$-adic expansion of the weight of~$B$. By \cite[Corollary~11.8]{OlssonCombi} we have
\begin{gather*}|\Irr_0(B)|=2^{\sum\limits_{i\ge 1}ia_{i-1}}.\end{gather*}

The cores of $2$-blocks are represented by the ``staircase'' partitions: $(k,k-1,\ldots,1)$ for $k\ge 0$. The Scopes classes are represented by those cores with $0\le k\le w-1$. In particular, there are just $w$ Scopes classes (unless $w=0$) and these cores are self-conjugate. Hence, for $\chi_\lambda\in\Irr(B)$ we have $\chi_{\lambda'}=\sgn\cdot\chi_\lambda\in\Irr(B)$ where $\sgn$ is the sign character and $\lambda'\in\MP(n)$ is conjugate to~$\lambda$. Since all $2$-regular elements lie in the alternating group, it follows that $[\chi_{\lambda'},\chi_\mu]^0=[\chi_\lambda,\chi_\mu]^0$. One can show further that $\bm{\lambda'}=((\bm{\lambda}_1)',(\bm{\lambda}_0)')$. Since $\chi_\lambda\in\Irr_0(B)$ implies that $\psi_{\bm{\lambda}}(1)$ is odd, we conclude that $|\bm{\lambda}_0|\ne|\bm{\lambda}_1|$ by~\eqref{deg} (unless $w=0$). In particular, $\chi_\lambda\ne\chi_{\lambda'}$ and it is enough to consider only half of the height $0$ characters. Let
\begin{gather}\label{M0}
\Delta:=\bigl\{\chi_\lambda\in\Irr_0(B)\colon |\bm{\lambda}_0|>|\bm{\lambda}_1|\bigr\}\qquad\text{and}\qquad M^0:=\bigl([\chi_\lambda,\chi_\mu]^0\bigr)_{\chi_\lambda,\chi_\mu\in\Delta}.
\end{gather}
Now for Morita equivalent $2$-blocks $B_1$ and $B_2$ with matrices~$M_1^0$ and~$M_2^0$ as above there exists a permutation matrix~$T^0$ such that $T^0M_1^0=M_2^0T^0$ (since~$M$ is invariant under the transpositions~$(\lambda,\lambda')$). If $w\le 28$ and $B_1$ and $B_2$ belong to different Scopes classes, we can show by computer that $T^0$ cannot exist. The computation for $w=29$ could not be completed due to memory restrictions.

Next, we further restrict ourselves to the case where $w$ is a $2$-power. Then $\chi_\lambda\in\Delta$ satisfies $\bm{\lambda}_1=()$ and $\bm{\lambda}_0$ is a hook partition by \eqref{deg} and the hook formula. Hence,
\begin{gather*}\Delta=\bigl\{\chi_\lambda\in\Irr(B)\colon \bm{\lambda}_0=\big(r,1^{w-r}\big),\ r=1,\ldots,w\bigr\}\end{gather*}
and $|\Delta|=w$. The matrix $X$ in Theorem~\ref{osima} becomes \begin{gather*}\bigl(\psi_{\bm{\lambda}}(g_{\bm{\mu}})\bigr)_{\chi_\lambda\in\Delta,\,\bm{\mu}\in\Gamma}=\bigl(\chi_{\bm{\lambda}_0}(s_{\bm{\mu}_1})\bigr),\end{gather*}
i.e., $X$ consists of rows of the character table of~$S_w$.

Now we determine the signs $\delta_p(\lambda)$ for $\chi_\lambda\in\Delta$. Let $\bm{\lambda}_0=\big(r,1^{w-r}\big)$, and let $\lambda_{(2)}=(k,k-1$, $\ldots,1)$ be the core of $B$. If $k$ is odd, then we add the partition
\begin{gather*}\big(2r,2^{\min\{k,w-r\}},1^{2\max\{0,w-r-k\}}\big)\end{gather*}
to $\lambda_{(2)}$ (componentwise) to obtain $\lambda$. In particular, $\delta_p(\lambda)=(-1)^{\max\{0,w-r-k\}}$. If, on the other hand, $k$ is even, then we add
\begin{gather*}\big(2(w-r+1),2^{\min\{k,r-1\}},1^{2\max\{0,r-1-k\}}\big)\end{gather*}
to $\lambda_{(2)}$ and conjugate afterwards to get $\lambda$. In this case we obtain $\delta_p(\lambda)=(-1)^{w-r+1+\min\{k,r-1\}}$. These observations made it possible to check Conjecture~\ref{main} for $w=32$. The next case $w=64$ is again out of reach.

In the final paragraph of this section we investigate further properties of $M$ which might be of interest for a potential theoretical proof.
We consider the diagonal of $M$ which is also the diagonal of
\begin{gather*}X\diag\big(2^{-l(\mu)}|{\C}_{S_w}(s_{\mu})|^{-1}\big)X^\text{t}.\end{gather*}
It is easy to see that $[\chi_\lambda,\chi_\lambda]^0=[\chi_\mu,\chi_\mu]^0$ if $\bm{\lambda}_0=(\bm{\mu}_0)'$. This means that the diagonal entries of $M$ come in pairs. On the other hand, our calculations indicate that $[\chi_\lambda,\chi_\lambda]^0>[\chi_\mu,\chi_\mu]^0$ whenever $\bm{\lambda}_0>\bm{\mu}_0>(w/2,1^{w/2})$ where $>$ denotes the lexicographical order. If this turns out to be true, then the permutation matrix~$T$ in Lemma~\ref{lem} must be a product of disjoint transpositions. In this case it suffices to find characters $\chi_\lambda,\chi_\mu\in\Irr(B_1)$ and $\chi_\sigma,\chi_\tau\in\Irr(B_2)$ such that $\bm{\lambda}_0=\bm{\sigma}_0=(\bm{\mu}_0)'=(\bm{\tau}_0)'$ and $\delta_p(\lambda)\delta_p(\mu)\ne\delta_p(\sigma)\delta_p(\tau)$. We do not know if this can be done in general.

\section{Implementation}\label{section6}

To check Conjecture~\ref{con1} for small primes, we have implemented the following steps in GAP by making use of the \texttt{hecke} package:
\begin{itemize}\itemsep=0pt
\item determine the cores of all Scopes classes for given $p$ and $w$ up to conjugation of partition,
\item use the Jantzen--Schaper formula to compute the decomposition matrices of Scopes classes (this is somewhat faster than Richard's formula),
\item compute the multisets of row sums and column sums for a given decomposition matrix,
\item partition the set of decompositions matrices according to their row sums and column sums (decomposition matrices in different parts of this partition cannot correspond to Morita equivalent blocks),
\item use the GAP command \texttt{TransformingPermutations} to distinguish decomposition matrices with the same row sums and column sums.
\end{itemize}

We have verified our code with the known decomposition matrices of $S_n$ with $n\le 18$. The applicability of the procedure above is only limited by the time it takes to run through all Scopes classes, but it can run in parallel. The computation of $M(11,2)=29{,}624$ took about one day on an Intel Xeon E5520 processor. We could not check Conjecture~\ref{con1} for $p=13$, since there are already $372{,}308$ classes to consider (hence $\approx10^{11}$ comparisons) and the decomposition matrix has size $104\times 90$ in each case. Similarly, computing $M(11,3)$ is hopeless, since this number is expected to be $4{,}209{,}504$.

On the other hand, our GAP programs for Conjecture~\ref{main} are mainly limited by the available physical memory. Our procedure here is given as follows:
\begin{itemize}\itemsep=0pt
\item compute the columns of the character table of $G_w$ corresponding to the elements $g_{\bm{\mu}}$ with $\bm{\mu}\in\Gamma$ using a code provided by Thomas Breuer (this uses much less memory than to compute the full character table via \texttt{CharacterTableWreathSymmetric}),
\item determine the cores of the Scopes classes as above,
\item compute the $p$-signs $\delta_p(\lambda)$ for $\chi_\lambda\in\Irr(B)$ with respect to the Scopes class just using the definition,
\item compute $M^0$ for each Scopes class,
\item compare the matrices $M^0$ using \texttt{TransformingPermutations},
\item for $p=2$ we only consider $\chi\in\Delta$.
\end{itemize}

All computations for Conjecture~\ref{main} take about five days on an Intel Xeon E5520 processor with 128GB memory.

Finally, for the exceptional case $(p,w)=(2,3)$ in Conjecture~\ref{main} we have computed the basic algebras of the blocks over $\FF_2$ in Magma and checked that these are isomorphic. All codes can be found on the author's \href{http://www.mathematik.uni-kl.de/agag/mitglieder/privatdozenten/dr-habil-benjamin-sambale/pub/}{homepage}.

\section{Concluding remarks}\label{section7}

In the following paragraph we indicate an alternative approach to the conjectures stated in the introduction. Let $T$ be the signed permutation matrix from Theorem~\ref{CR} with $TM_1=M_2T$. Suppose that $B_1$ and $B_2$ are Morita equivalent. Then we may assume that $M:=M_1=M_2$ by Lemma~\ref{lem}. Since $M$ is real and symmetric, $M$ and~$T$ are simultaneously diagonalizable. Since $M^2=M$ (see \eqref{form}), all eigenvalues of $M$ are $0$ and $1$. However, the situation gets more interesting if we restrict to characters of a~fixed height. As in Section~\ref{section5}, let $p=2$ and $w$ be a~$2$-power. Then the matrix $M^0$ in \eqref{M0} seems to have very interesting eigenvalues. For instance, the eigenvalues of $2^{12}M^0$ for $w=8$ are
\begin{align*}
2^{11},&& 2^{10},&& 2^9,&& 2^7\cdot3,&& 2^5\cdot3^2,&& 2^4\cdot3\cdot5,&& 2^3\cdot5^2,&& 5^2\cdot7.
\end{align*}
If we can show that the eigenvalues of $M^0$ are pairwise distinct, then $M$ and $T$ can be diagona\-li\-zed by the same real basis transformation. In this case $T$ has real eigenvalues and we conclude that $T$ is a product of disjoint transpositions (modulo signs).

For $p=3$ in Conjecture~\ref{main} it is not always enough to consider only height $0$ characters. Even worse for $p\ge 5$, Lemma~\ref{lem} cannot tell all Scopes classes apart. For example the $5$-blocks of weight $2$ with cores $\big(8,4^2,1^4\big)$ and $\big(9, 5^2, 2^3\big)$ have decomposition matrices
\begin{gather*}\begin{pmatrix}
1&\!.\!&\!.\!&\!.\!&\!.\!&\!.\!&\!.\!&\!.\!&\!.\!&\!.\!&\!.\!&\!.\!&\!.\!&.\\
.&\!1\!&\!.\!&\!.\!&\!.\!&\!.\!&\!.\!&\!.\!&\!.\!&\!.\!&\!.\!&\!.\!&\!.\!&.\\
1&\!1\!&\!1\!&\!.\!&\!.\!&\!.\!&\!.\!&\!.\!&\!.\!&\!.\!&\!.\!&\!.\!&\!.\!&.\\
.&\!.\!&\!1\!&\!1\!&\!.\!&\!.\!&\!.\!&\!.\!&\!.\!&\!.\!&\!.\!&\!.\!&\!.\!&.\\
.&\!.\!&\!.\!&\!1\!&\!1\!&\!.\!&\!.\!&\!.\!&\!.\!&\!.\!&\!.\!&\!.\!&\!.\!&.\\
.&\!1\!&\!1\!&\!.\!&\!.\!&\!1\!&\!.\!&\!.\!&\!.\!&\!.\!&\!.\!&\!.\!&\!.\!&.\\
1&\!.\!&\!1\!&\!.\!&\!.\!&\!.\!&\!1\!&\!.\!&\!.\!&\!.\!&\!.\!&\!.\!&\!.\!&.\\
.&\!.\!&\!1\!&\!1\!&\!.\!&\!1\!&\!1\!&\!1\!&\!.\!&\!.\!&\!.\!&\!.\!&\!.\!&.\\
.&\!.\!&\!.\!&\!1\!&\!1\!&\!.\!&\!.\!&\!1\!&\!1\!&\!.\!&\!.\!&\!.\!&\!.\!&.\\
.&\!.\!&\!.\!&\!.\!&\!.\!&\!.\!&\!1\!&\!1\!&\!.\!&\!1\!&\!.\!&\!.\!&\!.\!&.\\
.&\!.\!&\!.\!&\!.\!&\!.\!&\!1\!&\!.\!&\!1\!&\!1\!&\!1\!&\!1\!&\!.\!&\!.\!&.\\
.&\!.\!&\!.\!&\!.\!&\!.\!&\!1\!&\!.\!&\!.\!&\!.\!&\!.\!&\!1\!&\!1\!&\!.\!&.\\
.&\!.\!&\!.\!&\!.\!&\!.\!&\!.\!&\!.\!&\!.\!&\!1\!&\!1\!&\!1\!&\!1\!&\!1\!&.\\
.&\!.\!&\!.\!&\!.\!&\!.\!&\!.\!&\!.\!&\!.\!&\!.\!&\!1\!&\!.\!&\!.\!&\!1\!&1\\
.&\!.\!&\!.\!&\!.\!&\!.\!&\!.\!&\!.\!&\!.\!&\!.\!&\!.\!&\!.\!&\!1\!&\!.\!&.\\
.&\!.\!&\!.\!&\!.\!&\!.\!&\!.\!&\!.\!&\!.\!&\!.\!&\!.\!&\!.\!&\!.\!&\!.\!&1\\
.&\!.\!&\!.\!&\!.\!&\!.\!&\!.\!&\!.\!&\!.\!&\!.\!&\!.\!&\!.\!&\!1\!&\!1\!&1\\
.&\!.\!&\!.\!&\!.\!&\!.\!&\!.\!&\!.\!&\!.\!&\!1\!&\!.\!&\!.\!&\!.\!&\!1\!&.\\
.&\!.\!&\!.\!&\!.\!&\!1\!&\!.\!&\!.\!&\!.\!&\!1\!&\!.\!&\!.\!&\!.\!&\!.\!&.\\
.&\!.\!&\!.\!&\!.\!&\!1\!&\!.\!&\!.\!&\!.\!&\!.\!&\!.\!&\!.\!&\!.\!&\!.\!&.
\end{pmatrix}\!, \qquad
\begin{pmatrix}
1&\!.\!&\!.\!&\!.\!&\!.\!&\!.\!&\!.\!&\!.\!&\!.\!&\!.\!&\!.\!&\!.\!&\!.\!&.\\
.&\!1\!&\!.\!&\!.\!&\!.\!&\!.\!&\!.\!&\!.\!&\!.\!&\!.\!&\!.\!&\!.\!&\!.\!&.\\
1&\!1\!&\!1\!&\!.\!&\!.\!&\!.\!&\!.\!&\!.\!&\!.\!&\!.\!&\!.\!&\!.\!&\!.\!&.\\
.&\!.\!&\!.\!&\!1\!&\!.\!&\!.\!&\!.\!&\!.\!&\!.\!&\!.\!&\!.\!&\!.\!&\!.\!&.\\
.&\!.\!&\!.\!&\!1\!&\!1\!&\!.\!&\!.\!&\!.\!&\!.\!&\!.\!&\!.\!&\!.\!&\!.\!&.\\
.&\!1\!&\!1\!&\!.\!&\!.\!&\!1\!&\!.\!&\!.\!&\!.\!&\!.\!&\!.\!&\!.\!&\!.\!&.\\
1&\!.\!&\!1\!&\!.\!&\!.\!&\!.\!&\!1\!&\!.\!&\!.\!&\!.\!&\!.\!&\!.\!&\!.\!&.\\
.&\!.\!&\!.\!&\!1\!&\!.\!&\!.\!&\!.\!&\!1\!&\!.\!&\!.\!&\!.\!&\!.\!&\!.\!&.\\
.&\!.\!&\!.\!&\!1\!&\!1\!&\!.\!&\!.\!&\!1\!&\!1\!&\!.\!&\!.\!&\!.\!&\!.\!&.\\
.&\!.\!&\!.\!&\!.\!&\!.\!&\!.\!&\!.\!&\!.\!&\!.\!&\!1\!&\!.\!&\!.\!&\!.\!&.\\
.&\!.\!&\!.\!&\!.\!&\!.\!&\!.\!&\!.\!&\!.\!&\!.\!&\!1\!&\!1\!&\!.\!&\!.\!&.\\
.&\!.\!&\!.\!&\!.\!&\!.\!&\!.\!&\!.\!&\!.\!&\!.\!&\!.\!&\!1\!&\!1\!&\!.\!&.\\
.&\!.\!&\!.\!&\!.\!&\!.\!&\!.\!&\!.\!&\!.\!&\!.\!&\!1\!&\!1\!&\!1\!&\!1\!&.\\
.&\!.\!&\!.\!&\!.\!&\!.\!&\!.\!&\!.\!&\!.\!&\!.\!&\!1\!&\!.\!&\!.\!&\!1\!&1\\
.&\!.\!&\!.\!&\!.\!&\!.\!&\!1\!&\!.\!&\!.\!&\!1\!&\!.\!&\!.\!&\!1\!&\!.\!&.\\
.&\!.\!&\!.\!&\!.\!&\!.\!&\!.\!&\!1\!&\!.\!&\!1\!&\!.\!&\!.\!&\!.\!&\!.\!&1\\
.&\!.\!&\!.\!&\!.\!&\!.\!&\!.\!&\!.\!&\!1\!&\!1\!&\!.\!&\!.\!&\!1\!&\!1\!&1\\
.&\!.\!&\!.\!&\!.\!&\!.\!&\!.\!&\!.\!&\!1\!&\!.\!&\!.\!&\!.\!&\!.\!&\!1\!&.\\
.&\!.\!&\!1\!&\!.\!&\!1\!&\!1\!&\!1\!&\!.\!&\!1\!&\!.\!&\!.\!&\!.\!&\!.\!&.\\
.&\!.\!&\!1\!&\!.\!&\!1\!&\!.\!&\!.\!&\!.\!&\!.\!&\!.\!&\!.\!&\!.\!&\!.\!&.
\end{pmatrix}
\end{gather*}
respectively. Since the first matrix has five rows with only one non-zero entry and the second matrix has only four such rows, it is clear that the blocks cannot be Morita equivalent. However, both matrices yield the same matrix~$M$.

\subsection*{Acknowledgment}
The author likes to thank Susanne Danz for some valuable information on Scopes classes and Thomas Breuer for providing a GAP routine. Moreover, the author appreciates useful comments by an anonymous referee. Parts of this work were written while the author was in residence at the Mathematical Sciences Research Institute in Berkeley (Spring 2018) with the kind support by the National Science Foundation (grant DMS-1440140). This work is also supported by the German Research Foundation (\mbox{SA~2864/1-1} and \mbox{SA~2864/3-1}).

\pdfbookmark[1]{References}{ref}
\LastPageEnding

\end{document}